\documentclass{amsart}
\usepackage{amsmath,amsfonts,amssymb,amsthm}
\usepackage{mathtools}
\usepackage{graphicx} % Required for inserting images
\usepackage[english]{babel}
\usepackage[utf8]{inputenc}
\usepackage{hyperref}
\usepackage{csquotes}
\usepackage{xcolor}
\usepackage{pgfplots, tikz}
\usepackage[mathscr]{euscript}
\usepackage{graphicx}
\usepackage{subcaption}

\newtheorem{lemma}{Lemma}[section]
\newtheorem{corollary}[lemma]{Corollary}

\newtheorem*{maintheorem}{Main Theorem}
\newtheorem*{Corollary}{Corollary}

\theoremstyle{definition}
\newtheorem{definition}[lemma]{Definition}

\theoremstyle{remark}

\newcommand{\C}{\mathbb{C}}

\newcommand{\R}{\mathbb{R}}

\renewcommand{\Im}{\operatorname{Im}}
\renewcommand{\Re}{\operatorname{Re}}

\newcommand{\WDK}{\text{WDK}}
\newcommand{\EA}{\text{EA}}

\newcommand{\hide}[1]{}

\newcommand{\Newpage}{}

\begin{document}

\title{The Ehrlich--Aberth root finder is not generally convergent}
\author{Torben Grabbel, Nedim Srkalovic, Dierk Schleicher}

%\subjclass[2000]{65H04; 37F10; 37F80}

%\date{\today}

\maketitle

\begin{abstract}
    We show that the well-known Ehrlich--Aberth root finding method is not generally convergent. More precisely, we show that there is an open subset of polynomials of any degree $d\ge 7$, and an open subset of initial conditions, such that the iteration of the Ehrlich--Aberth method converges to an attracting periodic cycle of period $4$, rather than to a vector of roots. The proof is constructive in the sense that we present explicit values for the polynomial and the $4$-cycle.
    
    The search is highly non-trivial: we found no attracting cycles for lower degrees/periods, and the size of the attracting domain relative to a ``unit cube'' in the space of polynomials and starting vectors is on the order of $10^{-127}$.
\end{abstract}

\tableofcontents

\section{Introduction}

Finding roots of polynomials is an active field of research in mathematics, physics, engineering, and other disciplines, and it has important applications in all these areas. It has arguably several millennia of history, but it is still very active, and new methods continue to be developed.

It has been known since Gauss that every complex polynomial in one variable factors into roots, and since Ruffini--Abel that there is no solution formula for polynomials of degrees 5 or higher in terms of radicals. Therefore, most root finding methods are based on some sort of iteration.

Two main variants are root finders that iterate in $\C$ and converge to a single root (when they do), and those that iterate in $\C^d$ (where $d$ is the degree of the polynomial to be solved) and are supposed to converge to a $d$-vector of a permutation of all the roots.

Both variants have different advantages. The best known root finder in one variable is Newton's method, given by $N_p(z)=z-p(z)/p'(z)$ where $p$ is the polynomial to be solved. It is easy and usually fast to find a single root, but it is a challenge to find all $d$ roots. It has been known for a long time that it is not generally convergent (we define this term below): even the simple polynomial  $z^3-2z+2$ has an attracting cycle of period $2$ that attracts an open set of initial points. Nonetheless, polynomials of very large degrees have been factorized completely by algorithms based on Newton's method \cite{NewtonRobin2,TeraPolynomial}. Moreover, there is meanwhile a satisfactory theory both for where to start the iteration in order to find all roots \cite{NewtonInventiones,NewtonProbabilistic}, and how many iterations are required \cite{NewtonIterations,BAS,NewtonEfficient}. 

While it has been known for a long time that Newton's method is not generally convergent, it was shown by McMullen~\cite{McMullen} that there is no (holomorphic) iterative method in one complex variable that is generally convergent for polynomials of degree $4$ or greater.

Much attention is thus given to higher-dimensional root finders that iterate in $\C^d$ and converge to a $d$-vector that consists of a permutation of all the roots. Probably the two best known methods are those of Weierstrass (also known as Durand--Kerner method) and of Ehrlich--Aberth (see Section~\ref{Sub:BasicDefinitions} for precise definitions). The fact that these methods have multiple names is due to the fact that they have been re-discovered independently several times. We describe both methods below.

The Weierstrass--Durand--Kerner method (\WDK\ in short) is the oldest of these; it was introduced in 1891 \cite{WeierstrassOriginal} in the context of providing a constructive proof of the Fundamental Theorem of Algebra. The Ehrlich--Aberth method (\EA) is somewhat more complicated, but also known to be more efficient. It is about half as old as \WDK\ \cite{Ehrlich,Aberth}, and it is the basis for the professional root-finding implementation \textsc{MPSolve} by Bini and Robol \cite{MPSolve0,MPSolve}. 

The \WDK\ and the \EA\ methods have passed the test of time very well. They have been used in innumerable applications successfully, and they have always succeeded except in obvious cases of symmetry: if the given polynomial is real but has non-real roots, and the starting vector is real, then the iteration will always remain real and thus cannot converge to a root vector. This problem is not stable, so small perturbations will usually converge to a root vector. From a practical point of view, these methods are thus known to ``always converge''.

However, it was discovered a few years ago that the \WDK\ root finder is not generally convergent: for every degree $d\ge 3$, there is an open set of polynomials and an open set of initial vectors that converge to cycles of length $4$ or higher. This was shown in \cite{BernhardMichaelWeierstrass} based on groundbreaking work of B.~Reinke; we will describe some of the main ideas below (Section~\ref{Sub:Comparision_EA_WDK}).

Finding attracting cycles for the \EA\ method is a bigger challenge, for several reasons: the method itself is more complicated; and for every degree of freedom in the true search space there are two more eigenvalues to check (for two independent reasons) than for \WDK; see Section~\ref{Sub:Comparision_EA_WDK} for details. This makes the analysis much more complicated.

Nonetheless, here is our main theorem.
 
\begin{maintheorem}
	The Ehrlich--Aberth root finder has an attracting periodic orbit of period $4$ for an open set of polynomials of each fixed degree $d\ge 7$.
\end{maintheorem}

We present a specific polynomial and specific vector on an attracting cycle in Section~\ref{Sub:AttractingCycle}.

Of greatest structural relevance is the following corollary.

\begin{Corollary}
	The Ehrlich--Aberth root finder is not generally convergent.
\end{Corollary}

This result is the outcome of a six-year effort of our team to understand the global dynamics of the Ehrlich--Aberth method and the structure of its parameter space systematically in order to find (or rule out) the existence of attracting cycles. In this paper, we present the main result. For special cases we can rule out the existence of attracting cycles. These results will be published separately.

We should point out that the search is quite non-trivial, and a random search is not likely to find any attracting cycle. Indeed, we have not been able to find attracting cycles for lower degrees and periods than $(7,4)$, while for this combination of degree and period there might be up to $10^{13}$ periodic cycles of the type we consider, and for most polynomials and most starting vectors there is no attracting cycle.

A rough estimate (explained in Section~\ref{Sub:Stability}) shows that for the polynomial that features our attracting cycle, the diameter of its basin is, depending on directions, between $10^{-6}$ and $2\cdot 10^{-5}$. Similarly, within the space of degree 7 polynomials the diameters of how much the roots can be moved while preserving an attracting cycle are between $2\cdot 10^{-7}$ and $2\cdot 10^{-6}$. Overall, within the product space of polynomials and starting vectors, the locus of attracting dynamics occupies relative volume (compared to a unit cube) of less than $10^{-127}$, so we were looking for a tiny needle in a gigantic haystack. This explains why, despite heavy use in practice, an attracting cycle has never been found so far: even within a unit cube, a random polynomial and a random starting vector have (very roughly!) probability of less than $10^{-127}$ to converge to the attracting cycle we found; and the search space is certainly much larger than a unit cube.  (For comparison: the number of particles in the visible universe is estimated to be on the order of $10^{80}$.)

For completeness, we should mention that there are also efficient root finding schemes that do not depend on iteration, notably by square subdivision (as investigated in particular by Victor Pan~\cite{PanSquareSearch}), or homotopy methods.

\bigskip

\noindent
\textbf{Acknowledgments}. 
Over the several years that we have worked on this project, we have had most pleasant and useful discussions with a number of friends and colleagues, including Thomas Baldini, Dario Bini, Mika Grauel, Kilian Keen, Caspar Kiehn, Nicolae Mihalache, Victor Pan, Bernhard Reinke, Leonardo Robol, Oscar Scherz, and others. We would like to express our gratitude to all of them.

\medskip
\noindent
\textbf{Interaction with Artificial Intelligence}. 
We had worked in our team on finding attracting cycles for the Ehrlich--Aberth method for quite a few years, by various methods, ever since we discovered an attracting cycle for the Weierstrass method in the spring of 2020. When new tools based on AI developed, we started taking advantage of them in the late summer of 2026. The attracting cycle was thus found in close interaction with the AI system \emph{Claude}. This is thus an instance when AI helped experienced researchers make progress that would otherwise have been difficult to accomplish, rather than an AI system taking over. We have learned in the process that Claude had incorporated quite a collection of algorithms in the search, some with more success than others. --- It feels that we got access to a big machine, but without detailed operation manual on how to interact best, so this was an interesting learning experience for us. 

We would like to point out that the focus of this publication is on the dynamics of a well known root finding algorithm, not details on interaction with AI, or on details of implemented search algorithms. We will therefore only give a cursory overview on the latter aspects. See also  Section~\ref{Sub:RemarksComputers} on the use of computer clusters.

\Newpage

\section{Background}

\subsection{Basic definitions}
\label{Sub:BasicDefinitions}

To make the statements of our main theorem and its corollary precise we need a few preparations. Following McMullen's definition for one-dimensional root finders \cite{McMullen}, we make the following definition.

%more about best case scenario, i.e. open dense set with measure zero complement?

\begin{definition}[Convergence and general convergence]
    Let $p$ be a degree $d$ polynomial and $R:\C^d\dashrightarrow\C^d$ be a rational map. We say that $R$ is \emph{convergent} if $R^n(\vec{v})$ converges to some permutation of the roots of $p$ for an open dense subset of vectors $\vec{v}\in\C^d$.

    Let $\C_d[z]$ denote the set of degree $d$ polynomials and $\text{Rat}_k$ the set of degree $k$ rational functions in $d$ variables. An iterative root finder is a rational map $T\colon\C_d[z]\to\text{Rat}_k$ from degree $d$ polynomials to degree $k$ rational maps. It is called \emph{generally convergent} if there exists an open dense set $U\subseteq\C_d[z]$ such that for $p\in U$ the rational map $T(p)$ is convergent. 
\end{definition}

We note that our Main Theorem disproves general convergence in a very precise sense, by providing an open set of polynomials where the \EA\ root finder does not converge, but has a different kind of simple dynamics. 

Now we define the \WDK\ and \EA\ root finders.

\begin{definition}[Weierstrass--Durand--Kerner method (WDK)]
    Let $p$ be a monic degree $d$ polynomial. The \WDK\ iteration step is defined as: \[
        \WDK_p:\C^d\dashrightarrow\C^d; v_i\mathrel{\mapstochar\dashrightarrow} v_i-\frac{p(v_i)}{\prod_{j\neq i}(v_i-v_j)}
        \;.
    \]
        
\end{definition}

The idea is to pretend that, for each root separately, all the other roots have already been found; then the quotient equals $v_i$ minus the missing root (see \cite[Section~3]{BernhardMichaelWeierstrass} for a more detailed description). 

\begin{definition}[Ehrlich--Aberth method (\EA)]\label{def:EA_method}
    Let $p$ be a (monic) degree $d$ polynomial. The Ehrlich--Aberth iteration step is defined as: \[\EA_p:\C^d\dashrightarrow\C^d;\ v_i\mathrel{\mapstochar\dashrightarrow} v_i-\frac{\frac{p(v_i)}{p'(v_i)}}{1-\frac{p(v_i)}{p'(v_i)}\cdot\sum_{k\neq i}(\frac{1}{v_i-v_k})}\;.
    \]
\end{definition}

The motivation is similar as for \WDK, except that in each coordinate $v_i$ a single one-dimensional Newton step is applied. More on the interpretation of both methods can be found for instance in \cite[Sections~3 and 4]{Reinke}. 

In particular, one can view the Ehrlich--Aberth method as $d$ Newton orbits in parallel that are synchronized in such a way that when one orbit converges to a simple root, other orbits are pushed away so they search different roots.

We should clarify that in our approach, all $d$ coordinates are updated simultaneously (the \emph{Jacobi variant}). This is in contrast to the \emph{Gauss--Seidel variant} where the coordinates are updated cyclically (that is the variant used in MPSolve, but it is not as easy to investigate as a dynamical system).

Our approach is guided by viewing the Ehrlich--Aberth iteration $\EA_p$ for the polynomial $p$ as a dynamical system defined on a subset of $\C^d$. The entire method $\EA_{\cdot}$ can be thought of as a space of dynamical systems parametrized by all polynomials $p\in\C_d[z]$ resp.\ by their roots. 

\subsection{Symmetries}

An obvious symmetry of the system is that permutations of the vector entries commute with iteration of \EA. There is a somewhat more subtle symmetry, as follows.

\begin{lemma}[Möbius equivariance, \protect{\cite[Lemma~3.2]{Reinke}}]
    Let $p(z)=\prod_{i=1}^d(z-z_i)$ be a polynomial and $M$ be a Möbius transformation so that all $M(z_i)\neq\infty$. Define $\tilde{p}(z)=\prod_{i=1}^d(z-M(z_i))$. Then the Ehrlich--Aberth method satisfies the following equivariance property: \[
    M(\EA_p(v_1,\dots,v_d))=\EA_{\tilde{p}}(M(v_1),\dots,M(v_d)).
    \]
\end{lemma}

In particular, $M$ conjugates the action of $\EA$ for $p$ to that of $\tilde p$, and thus sends all periodic orbits of $\EA_p$ to periodic orbits of $\EA_{\tilde p}$ with equal period and equal eigenvalues.

This property allows us to place three roots of the polynomial $p$ freely, which reduces the complexity of the parameter space significantly; see the discussion in Section~\ref{Sub:SearchSpace}. 

\subsection{Attracting cycles and eigenvalues}

A possible obstruction for a root finder to be generally convergent is to possess attracting cycles (this is not the only kind of conceivable obstruction, but it seems to be the easiest one to find, and the only one known to exist for \WDK\ and \EA). A cycle of period $n$ for the Ehrlich--Aberth method for a polynomial $p$ is a set of vectors $\{\vec{v}^1,\dots,\vec{v}^n\}$ satisfying $\EA_p(\vec{v}^i)=\vec{v}^{i+1}$ for $i\in\{1,2,\dots,n\}$ (indices modulo $n$). A cycle is 
\emph{attracting} if there exists a neighborhood $U$ containing $\{\vec{v}^1,\dots,\vec{v}^n\}$ with $\EA_p^{i+1}(\overline{U})\subset \EA_p^i(U)$ and $\bigcap_{i=0}^\infty \EA_p^i(U)=\{\vec{v}^1,\dots,\vec{v}^n\}$. The discrete time version of the Hartman--Grobman theorem (see for instance \cite[Section~1.4]{GuckenheimerHolmes}) states that the linearization of the first return map gives the following sufficient criterion. 

\begin{lemma}[Eigenvalues and attracting cycles]
    A cycle of period $n$ is attracting if the eigenvalues $\rho_1,\dots,\rho_d$ of the Jacobian of the first return map $D(\EA_p^n)(v_i)$ for some $i\in\{1,\dots,n\}$ all have  absolute values strictly less than $1$.
\end{lemma}

All points $v_i$ on the periodic cycle have conjugate Jacobians, so the criterion is satisfied for all points on the same cycle simultaneously. 

For a periodic cycle $\{v_1,\dots,v_n\}$, we denote by $\rho=\rho(\{v_1,\dots,v_n\})$ the absolute value of the largest eigenvalue(s). Our goal is to find a cycle with $\rho<1$.

\subsection{Multiple leading eigenvalues}
\label{Sub:MultipleLeadingEigenvalues}
It turns out that, at least generically, one can expect that for a polynomial that locally minimizes the leading eigenvalue, at least $d-2$ eigenvalues have equal absolute values. This is because the parameter space is ($d-3$)-dimensional, and when less than $d-2$ eigenvalues have equal minimal absolute values, one can locally minimize the leading eigenvalue. (The simplest example is when there is a single degree of freedom: if the leading eigenvalue is simple, then one should be able to optimize it until it hits at least one more eigenvalue with equal absolute value.)

Our attracting cycle presented below actually has $d-1$  eigenvalues with equal absolute values, one more than expected, while the last eigenvalue is rather close also.

\subsection{Cycles as algebraic varieties}\label{sec:cyc_to_var}

Fix a degree $d$ and a cycle period $n$ and introduce coordinates $r_1,\dots,r_d$ for the roots and $v_i^j$ for $i\in\{1, 2, \dots, d\}$ and $j\in\{1, 2, \dots, n\}$ for the entries of the cycles. The cycle condition $\vec{v}\,^{j+1}=\EA_p(\vec{v}\,^j)$ (coordinate wise) implies a system of rational equations. The intersection of the solution sets is the set of cycles (admitting solutions that come from cycles of lower periods, which exist but are relatively few). By clearing the denominators in the rational equations one obtains a system of polynomials whose corresponding variety contains the set of cycles. In general, this variety will also contain degenerate solutions corresponding to two entries being equal, entries being equal to roots of polynomials, or a denominator in the definition of \EA\ being equal to zero. Thus we can describe the degenerate solutions using their own polynomial equations, and remove them. For fixed values of $r_1,\dots,r_d$ there are generically only a finite number of cycles by Bézout's theorem. 

The equations described here are the ones that we try to solve numerically. 

\subsection{Comparison between \EA\ and \WDK, and the approach in \cite{BernhardMichaelWeierstrass}}
\label{Sub:Comparision_EA_WDK}
A priori, the \WDK\ and the \EA\ methods are similar: they converge fast to a vector of roots in a certain neighborhood of the root vectors; for simple roots, the convergence is quadratic for \WDK\ and cubic for \EA.

The \WDK\ method has the advantage that the dynamics always takes place in a co-dimension $1$ hyperplane, so among the $d$ eigenvalues, one is always zero. This simplifies the search for attracting cycles, but this is not so for \EA.

A second relevant difference between \WDK\ and \EA\ is that the former only respects affine M\"obius transformations, hence only a $2$-dimensional subspace of polynomials yields conjugate dynamics. Therefore, the smallest non-trivial parameter space for \WDK\ is formed by cubic polynomials, while for \EA\ it is formed by quartic polynomials. Together with the fact that \EA\ in itself is more complicated, this yields three relevant reasons why the search for attracting cycles is more challenging for the \EA\ method.

The approach in \cite{BernhardMichaelWeierstrass} was essentially algebraic. It was sufficient to consider cubic polynomials. For $2$-cycles and $3$-cycles, it was possible using computer algebra systems to calculate explicit algebraic relations between the two non-zero eigenvalues, and to conclude that the two eigenvalues can never be simultaneously smaller than $1$ in absolute value. General $4$-cycles are already far out of reach for the currently available computer algebra systems, but a special type of permutation $4$-cycles could still be handled (as in Section~\ref{Sub:SymmetrySubspaces}), and it was this kind of cycles that yielded an attracting cycle. 

For the \EA\ method, the simplest non-trivial case are period $2$ cycles for degree $4$ polynomials, and these can be only partially handled by computer algebra systems. Beyond this, we had to resort to more numerical and even probabilistic methods.

\Newpage

\section{The main steps for the search}

\subsection{The search space} 
\label{Sub:SearchSpace}
We are searching for a \emph{polynomial} that has a \emph{periodic cycle} that is \emph{attracting}. The main combinatorial data are thus
\begin{enumerate}
    \item the degree $d$ of the polynomial $p$,
    \item the period $n$ of a periodic cycle (and possibly a symmetry type of the cycle), and
    \item a choice among all cycles of this period.
\end{enumerate}

It suffices to consider monic polynomials. Moreover, it suffices to consider the case that all roots are distinct and hence simple: if a polynomial with multiple roots has an attracting cycle, then this cycle is stable so that small perturbations with all simple roots still have an attracting cycle of the same period. We will thus assume that all roots are simple. 

Monic polynomials are defined by their (unordered) sets of roots. Since any two M\"obius equivalent sets of roots lead to conformally conjugate \EA\ maps, one can fix three roots of the polynomials at arbitrary complex positions, for instance at $0$, $1$, and $-1$. Equivalently, one can require that there is an $r\not\in\{0,1,-1\}$ such that four roots are at $1$, $-1$, $r$, and $-r$. 

In any case, the space of \EA\ maps of degree $d$ polynomials is complex ($d-3$)-dimensional. 

Since all cubic polynomials (with simple roots) lead to conformally equivalent \EA\ maps, the simplest infinite space to consider is the space of degree $d=4$ polynomials, and to search for $2$-cycles. The number of periodic points for any pair $(d,n)$ is finite, but it may be huge, and one would expect it to increase at least exponentially both in $d$ and $n$.
To determine good (upper or lower) bounds on this number seems to be difficult. Reinke and Silversmith~\cite{ReinkeSilversmith1,ReinkeSilversmith2} provide a formula for the number of periodic points of a generalization of the Ehrlich--Aberth method, which leads to an upper bound for the periodic orbits of \EA\ itself. These bounds are difficult to even evaluate, and it is not clear how much they over-estimate the actual counts. Here are a few values based on the formulas in \cite{ReinkeSilversmith1,ReinkeSilversmith2}, for various degrees $d$ and periods $n$; these numbers grow exponentially in $n$, and super-exponentially~in~$d$.
\[
\begin{array}{|r|rrr|}
\hline
	& n=2 & n=3  & n=4 \\ \hline
	d=4 \rule{0pt}{11pt}&6.0\cdot 10^4 & 8.2\cdot 10^7 & 5.4\cdot 10^{10}  \\
	d=5 & 2.3\cdot 10^7 & 5.9\cdot 10^{11} & 6.8\cdot 10^{15} \\
	d=6 & 1.3\cdot 10^{10} & 7.1\cdot 10^{15} & 1.8\cdot 10^{21} \\
	d=7 & 9.1\cdot 10^{12} & 1.3\cdot 10^{20} & 8.6\cdot 10^{26} \\ \hline
\end{array}
\]

Many periodic orbits are equivalent: by cyclic renumbering or by permuting coordinates. In general, the number of equivalent orbits is thus $d!\cdot n$, so the number of equivalence classes is correspondingly smaller.

In any case, when the search for periodic orbits for any pair $(d,n)$ has not led to attracting cycles, one can either increase the degree or the period (or both).

The methods of Reinke--Silversmith do not currently lead to estimates for cycles of permutation type. Specifically for $4$-cycles of permutation type, we suspect that the bounds are similar to those of general $2$-cycles: both are described by equations of the form $\EA^2(v)=P v$, where $P$ is the identity for $2$-cycles and a certain permutation matrix for permutation $4$-cycles. Recall that the Reinke--Silversmith bound is exact for a generalization of \EA\, but just an upper bound (not very tight) for \EA\ itself. Specifically, we thus use $9.1\cdot 10^{12}$ as a heuristic upper bound for the number of permutation $4$-cycles for degree $7$.  

When the search for a particular degree and period has not led to an attracting cycle, one can increase either of the two variables. Increasing the degree also increases the number of eigenvalues that must be in the unit disk, so as a first step it seems easier to increase the period instead. We have experimented with increasing either variable. 

\subsection{Symmetry subspaces}
\label{Sub:SymmetrySubspaces}

One strategy that has led to success for the Weierstrass method \cite{BernhardMichaelWeierstrass} is to restrict to periodic orbits that satisfy particular symmetries. A simple example of a symmetry cycle is a $2$-cycle $\{v_1,v_2\}$ such that $v_2=-v_1$, or such that $v_2$ is a permutation of $v_1$ of order $2$ without fixed elements, i.e.\ a pair exchange of coordinates. 

Such restrictions have the advantage that they provide intermediate conditions between increasing the search space by increasing the periods, while generating more manageable sets of solutions. 

In \cite{BernhardMichaelWeierstrass}, work was always on cubic polynomials ($d=3$), and it was shown that periodic orbits of period $2$ and $3$ are never attracting. Period $4$ was out of reach for the available methods, and the compromise was found to consider $4$-cycles with an additional symmetry after the second iterate. Interestingly, the attracting cycle found for \EA\ also has this type of symmetry.

\subsection{Local refinements}
Once a candidate cycle for a given polynomial has been identified (by any means, systematic or random search), such that the leading eigenvalue is ``not too far'' from $1$, then
one can try to reduce the leading eigenvalue by local optimization. A natural method is of course gradient descent. 
However, as described in Section~\ref{Sub:MultipleLeadingEigenvalues}, one expects minima for leading eigenvalues to be realized so that at least $d-2$ eigenvalues have coinciding absolute values. Therefore, at the points of interest in the space of polynomials, when leading eigenvalues are minimal, gradient methods fail to work, and one has to use higher order methods (not necessarily for finding attracting cycles, but when searching for smallest eigenvalues). (Note that, in a neighborhood of multiple eigenvalues, the behavior of the eigenvalues is typically root-like, so it is not even Lipschitz.)

We should mention that many of the candidate polynomials with small leading eigenvalues have not been found by refining initial polynomials for which the eigenvalues were already small, but often starting from where the eigenvalues were of moderate size.

\subsection{Increasing the degree}

Attracting cycles can be promoted from any degree to all higher degrees, in  analogy to \cite[Lemma 3.7]{BernhardMichaelWeierstrass}.

\begin{lemma}[Increasing degrees]
\label{Lem:IncreasingDegrees}
    Let $p$ be a monic polynomial of degree $d$ and let $\alpha\in\C$ be different from all the roots. Set $\tilde{p}(z)=(z-\alpha)p(z)$.
    \begin{enumerate}
        \item For a vector $\vec{v}=(v_1,\dots,v_d,\alpha)$ with pairwise distinct entries, the Jacobi matrix $D(\EA_{\tilde{p}})(\vec{v})$ has the form \[
        D(\EA_{\tilde{p}})(\vec{v})=\begin{pmatrix}
            D(\EA_{p})((v_1,\dots,v_d)) & *_{d\times1} \\
            0_{1\times d} & 0
        \end{pmatrix}.
        \]
        \item If $\vec{q}\in\C^d$ is a periodic point of $\EA_p$ of period $n$ such that all eigenvalues of $D(\EA_p^n)(\vec{q})$ have absolute values strictly less than $1$, and $\alpha$ is different from all entries in $\vec{q}$ and in its orbit, then $(\vec{q},\alpha)\in\C^{d+1}$ is a periodic point of $\EA_{\tilde{p}}$ of period $n$, such that all eigenvalues of $D(\EA_{\tilde{p}}^n)((\vec{q},\alpha))$ have absolute values strictly less than $1$.
    \end{enumerate}
\end{lemma}
\begin{proof}
    $(1)$ follows immediately from a computation.

    For $(2)$ apply the chain rule to compute the Jacobian $D(\EA^n_{\tilde{p}})((\vec{q},\alpha))$ and observe that it is of the form $\begin{pmatrix}
        D(\EA_p^n)(\vec{q}) & *_{d\times 1} \\
        0_{1\times d} & 0
    \end{pmatrix}$ where $*_{d\times1}$ is an unspecified vector in $\C^d$. Thus the eigenvalues of $D(\EA_{\tilde{p}}^n)((\vec{q},\alpha))$ are those of $D(\EA_{p}^n)(\vec{q})$ together with~$0$.
\end{proof}

\begin{corollary}[Increasing degrees improves eigenvalues]
\label{Cor:IncreasingDegrees}
	If there is a degree $d$ polynomial for which the \EA\ method has a periodic orbit of period $n$ and with leading eigenvalue with absolute value $\rho$, then there are also polynomials of all degrees greater than $d$ that have period $n$ orbits with leading eigenvalues with absolute value at most $\rho$. 
\end{corollary}

Lemma~\ref{Lem:IncreasingDegrees} and Corollary~\ref{Cor:IncreasingDegrees} remain true for periodic orbits of given symmetry type.

\subsection{The structure of the search}

There are at least three different ways to organize the search for attracting cycles (or to show non-existence) for a given degree, period, and possibly symmetry type:

\begin{enumerate}
	\item algebraic;
	\item numeric exhaustive; and
	\item probabilistic.
\end{enumerate}

All three methods come in different versions, and the borders between them are not necessarily well defined.

The most satisfactory of the three methods is the \emph{algebraic} one: describe the space of periodic points of given type by algebraic equations, discuss how the respective varieties decompose into components, calculate the eigenvalues as algebraic functions of the parameters (the polynomial roots). This is the method that has led to success for the \WDK\ method in \cite{BernhardMichaelWeierstrass}. For \EA, we have been able to do this successfully for many $2$-cycles of degree $4$ when they had suitable orbit symmetries. This has led, in various cases, to algebraic proofs of non-existence of attracting cycles of such types. This is an interesting project that we plan to describe separately; it has not led to the attracting cycle that we describe here.

The second method is \emph{numeric}: the periodic cycles of desired degree, period, and symmetry type are described by algebraic equations that are too complicated to be analyzed algebraically, but the (often very many) solutions can be computed numerically (if not all then many of them) and continued analytically throughout search space. Here the search space is limited in principle only by computing power, far beyond the limits of current computer algebra methods. When the search space is sufficiently manageable, one can use complex interval arithmetic to prove that, for certain domains in the search space, attracting cycles would have been found if they existed. For the ever-more huge search spaces described above, formal negative certificates have also become infeasible. 

This method was used for degree $4$ and period $2$ essentially up to saturation. 

For higher degrees and periods, a \emph{probabilistic} search was implemented. The probabilistic aspect refers to both the search for periodic cycles of given period and symmetry type for a fixed polynomial, and to the search through the space of polynomials. As described in Section~\ref{Sub:SearchSpace}, for any fixed polynomial and given period and cycle type, the number of cycles for a fixed polynomial may be huge, and very hard to estimate. It is unrealistic to find all periodic cycles systematically. Therefore, periodic orbits were searched for randomly, and it was tried to estimate how close the search was to saturation depending on how many separate orbits had been found exactly once, or exactly twice. 

The probabilistic search through the space of polynomials involved a Metro\-po\-lis-type algorithm: given a candidate periodic cycle, a fixed number of different nearby polynomials are explored (the $d-3$ roots are changed in different directions and distances, then a local optimization tries to continue analytically the previously discovered cycle, and then to find the best nearby polynomial where the leading eigenvalue is least). If this leads to a polynomial with better eigenvalue, this is the new base point and the step size is slightly increased. If not, then with a certain probability (depending on an ``external temperature'') the new candidate is still accepted. Otherwise, the previous candidate is continued to be explored with somewhat smaller step size. This strategy was followed for a fixed number of iterations, and then re-started at different random points. 

A probabilistic search was run for general $2$-cycles for degrees 5--10, but no attracting cycles have been found. 

Next the search went to $3$-cycles and ordinary $4$-cycles for degree $4$ polynomials; no attracting cycles were found. 

Finally, the search moved to transposition $4$-cycles: these are $4$-cycles where the second iterate is a permutation of the original vector. Therefore, the algebraic complexity is comparable to $2$-cycles, but there are additional candidates to investigate. 

For such cycles, degrees $4$ and higher were investigated. 
And indeed, for this cycle type the attracting orbit was found. 
Unlike the satisfactory result in \cite{BernhardMichaelWeierstrass} for \WDK\  (the combination $(3,4)$ for degree and period is minimal), we cannot make any claim on either the minimality of the degree or the period.

\subsection{The attracting $4$-cycle for degree $7$ and greater}
\label{Sub:AttractingCycle}

For the monic polynomial $p(z)=\prod_{i=1}^7(z-r_i)$ of degree $7$ with roots as in Table~\ref{Table:RootsDegreeSeven}, there is an attracting $4$-cycle near a starting vector that is also presented in the table, together with the seven eigenvalues.

\begin{figure}[htbp]
      \centering
      \framebox{
      \includegraphics[width=0.61\linewidth]{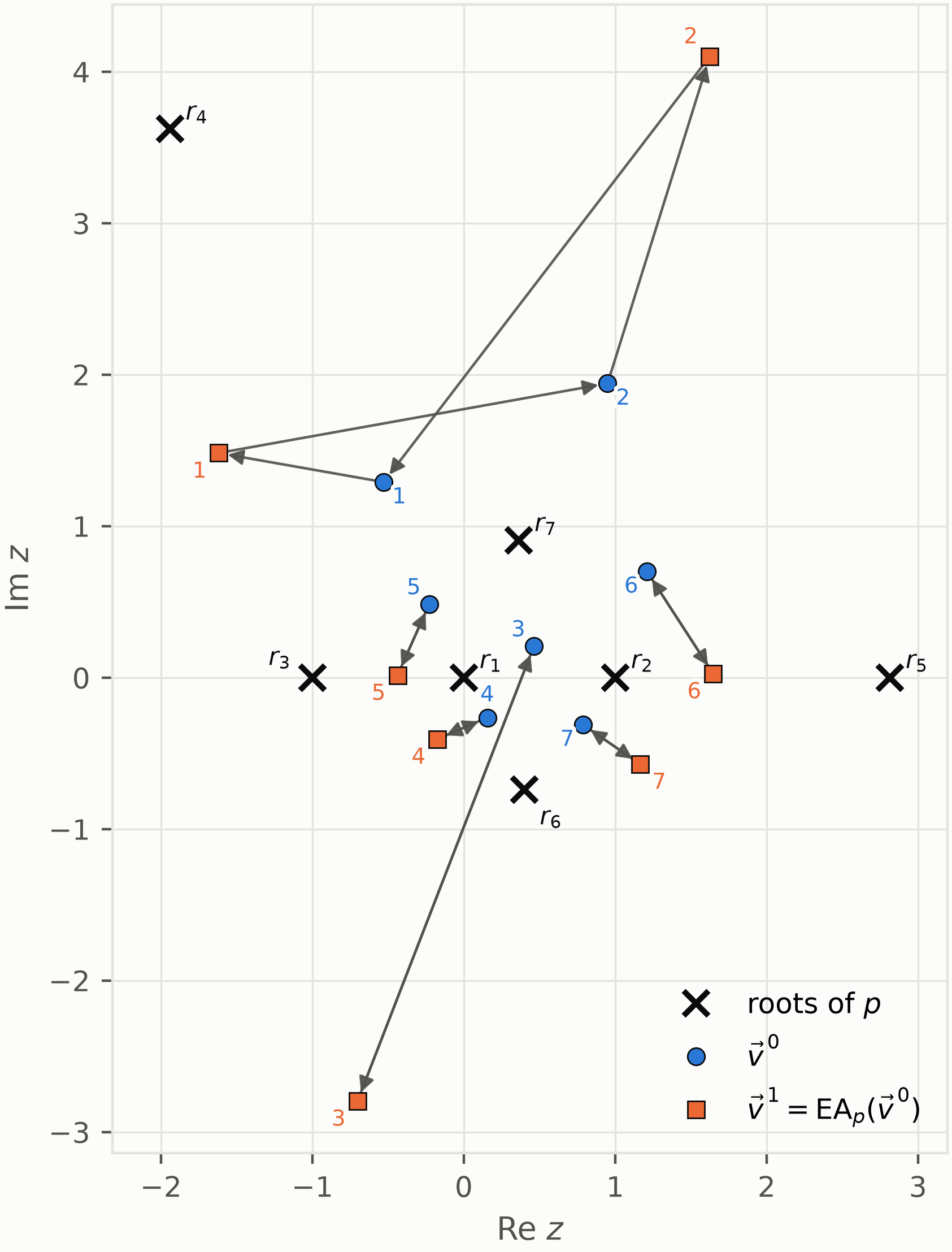}}
      \caption{The roots of $p$ and the two configurations $v^0$ (blue circles) and $v^1$ (red squares) of the attracting 4-cycle. The label at a point is the index of the approximation through the four steps $v^0 \mapsto v^1 \mapsto \tilde v^0 \mapsto \tilde v^1 \mapsto v^0$, where $\tilde v^0$ and $\tilde v^1$ differ from $v^0$ and $v^1$ only by a single pair exchange.}
      \label{fig:orbit_plane}
  \end{figure}

In Figure~\ref{fig:orbit_plane} we show, projected to a single complex plane, the positions of the seven roots, and of the attracting $4$-cycle. It turns out that five of the seven coordinates in this $4$-cycle actually have period $2$, while the remaining two coordinates switch positions after $2$ iterations: the attracting orbit thus has a very special symmetry type. 

This proves (numerically) the existence of an attracting $4$-cycle for a particular degree $7$ polynomial. We have  verified it numerically with complex interval (box) arithmetic using \cite{FLINT}. The precise starting values lie in boxes with side length $10^{-25}$ around those in Table~\ref{tab:values}. The absolute values of the eigenvalues are bounded by $0.7353184$. Moreover, we are grateful to Bernhard Reinke for an independent verification of this attracting cycle.

It is routine to upgrade the attracting cycle from degree~$7$ to all larger degrees: this follows from Corollary~\ref{Cor:IncreasingDegrees}.

\begin{table}[htbp]
    \[
    \begin{array}{r|rr} & \text{roots, real part} & \text{roots, imaginary part} \\ \hline
         1& 0 & 0  \\
 2& 1& 0  \\
 3& -1 & 0 \\
 4& -1.941586496092986502404187& 3.620929464488302329527869 \\
 5& 2.814011832575733951244956& -0.000008723882047762163862553685 \\
 6& 0.3976539063762122849787772 & - 0.7412693453568267587527885 \\
 7&  0.3627358400914195701858489   & 0.9041332040615440712727491 \\ \hline \hline
 & \text{starting vector, real part} \rule{0pt}{11pt} & \text{starting vector, imaginary part} \\ \hline
 1& -0.5290475728682283223117312& 1.2902311378633051304550690 \\
 2& 0.9496391591118217234183813& 1.9418382183527572037723937 \\
 3& 0.4628042056093357743829914& 0.2073255914500557057832052 \\
 4& 0.1551005076737854148786264& -0.2645761308773563410652116 \\
 5& -0.2276879008154907752369263& 0.4852564567782651106722482 \\
 6& 1.2088953169849776154606041& 0.7011279978927359339012465 \\
 7& 0.7858561209653443509429336& -0.3107770535938863118230804 \\ \hline\hline
 & \text{eigenvalue, real part} & \text{eigenvalue, imaginary part} \rule{0pt}{11pt} \\ \hline
 1& -0.061414018245771& 0.732749153879621 \\
 2& 0.734305074664565& -0.038566183526299 \\
 3& 0.731705822984688& 0.072785637293985 \\
 4& -0.013838166699863& -0.735186821972330 \\
 5& 0.734545660868135& -0.033670815529541 \\
 6& -0.106248772603267& -0.727600206760092 \\
 7& -0.060527920455661& 0.730400673305728 \\\hline
 \end{array}
\]
    \caption{The roots of the degree $7$ polynomial with an attracting $4$-cycle (top; exact values) and an element of its attracting $4$-cycle (center; 25 decimals, all digits certified). Bottom: the seven eigenvalues, all with absolute values about $0.7353$, except the last one is around $0.7329$.}
    \label{Table:RootsDegreeSeven}
    \label{tab:values}
\end{table}

\subsection{Stability and relative size}
\label{Sub:Stability}

Every attracting cycle is of course stable under perturbations: this orbit attracts a whole neighborhood of starting vectors in $\C^7$. Moreover, the cycle and its eigenvalues depend continuously on the parameters (the positions of the roots): therefore, for all sufficiently nearby polynomials there are attracting cycles of the same type. This proves, for every degree $d\ge 7$, the existence of an open set of polynomials, and an open set of starting vectors, that converge to an attracting $4$-cycle: so this concludes the proof of our Main Theorem, and thus 
disproves general convergence of the \EA\ method.

However, this phenomenon is ``exceedingly small'', in the following sense: rather small perturbations both of the positions of the roots, as well as the initial points on the iteration vector, lead to a loss of the attracting cycle. 

\begin{figure}[htbp]
\includegraphics[width=0.49\textwidth]{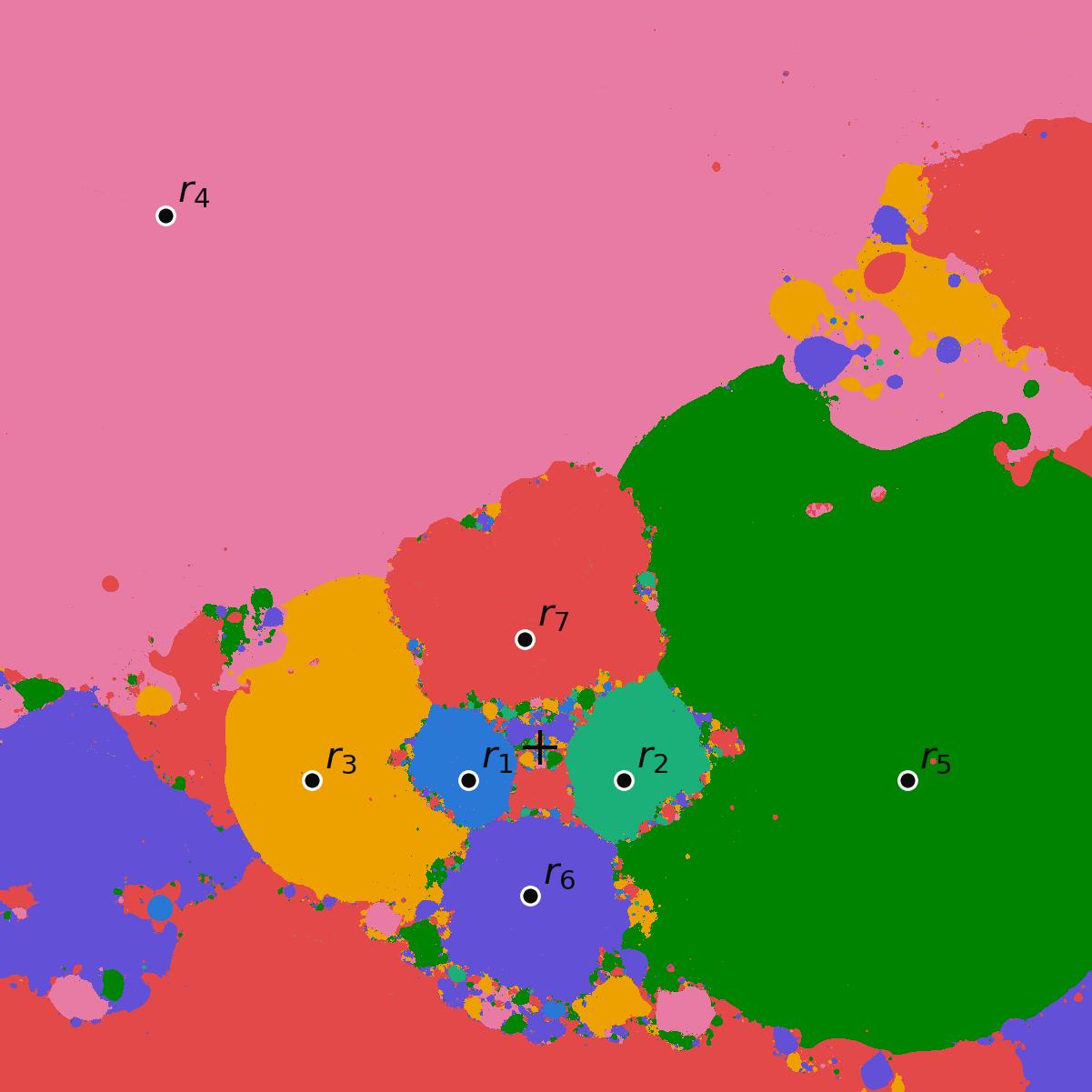}\hspace{2pt}
\includegraphics[width=0.49\textwidth]{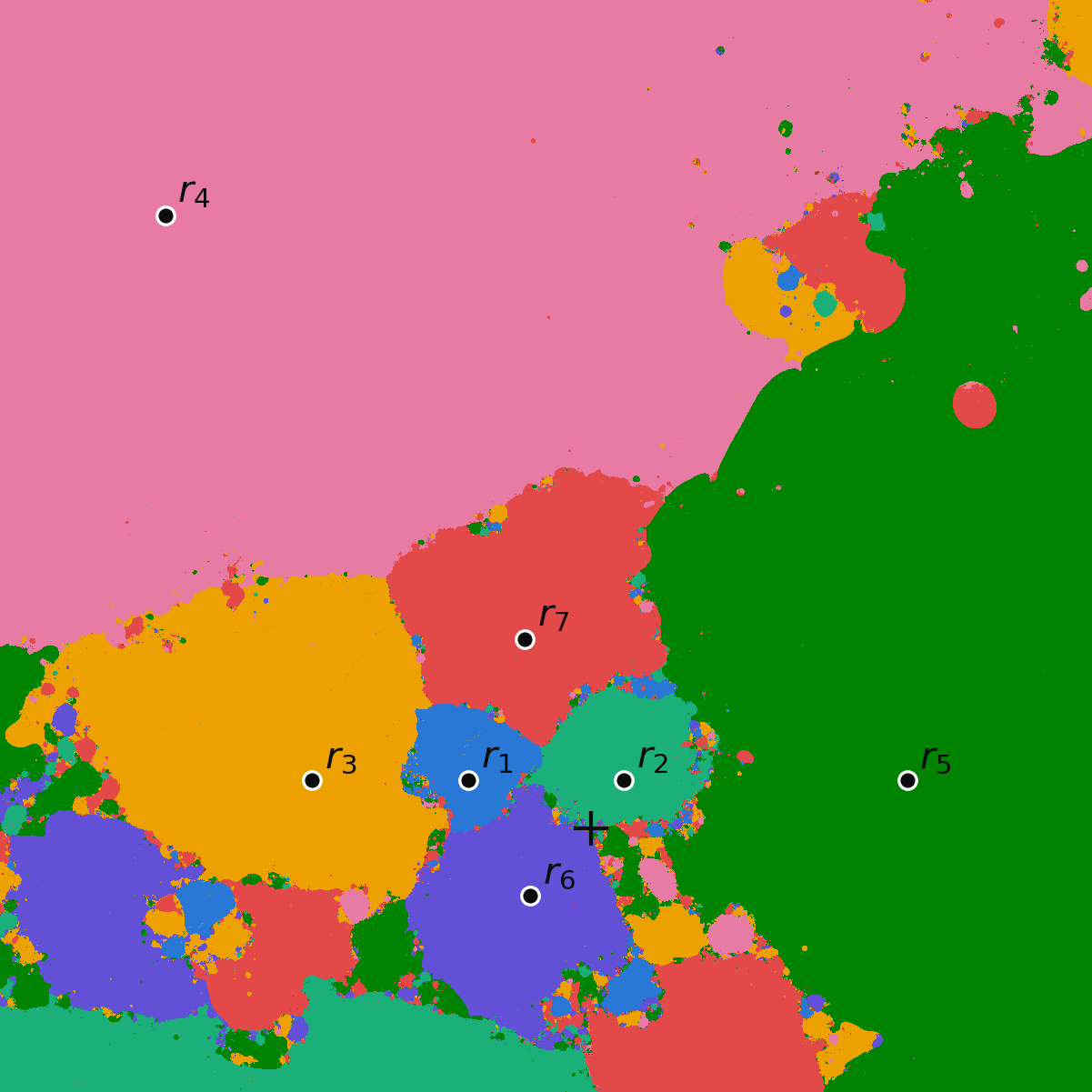}
	\caption{One-dimensional complex slices through the dynamics space. A single dynamical variable is changed in a neighborhood of the value corresponding to the attracting cycle; all other variables are held fixed. Different colors show the root that the coordinate under display converges to (the other coordinates may converge to the other roots in different permutations). The side length is $7$.}
	\label{Fig:DynamicSlices1}
\end{figure}

\begin{figure}[htbp]
\includegraphics[width=0.49\textwidth]{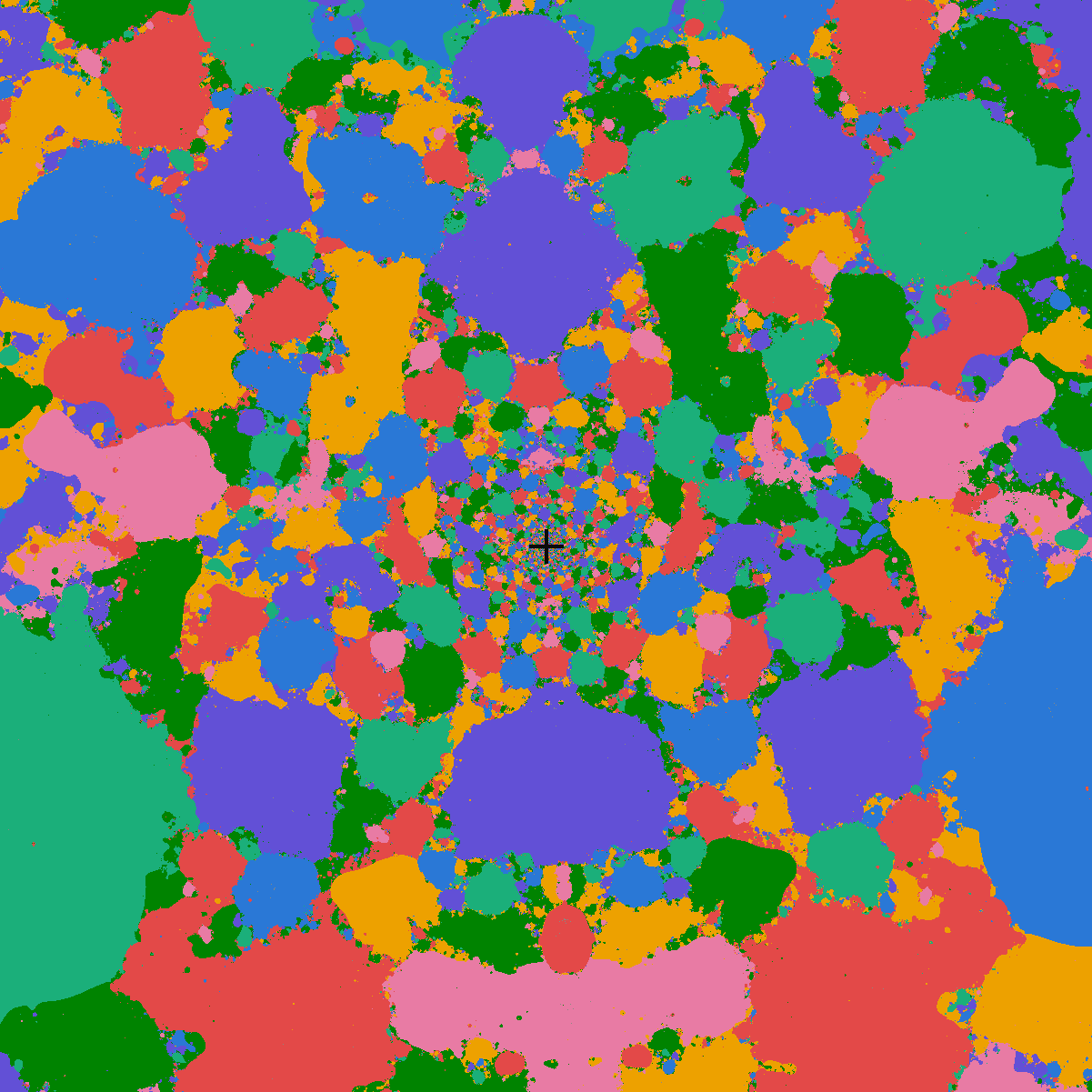}\hspace{2pt}%
\includegraphics[width=0.49\textwidth]{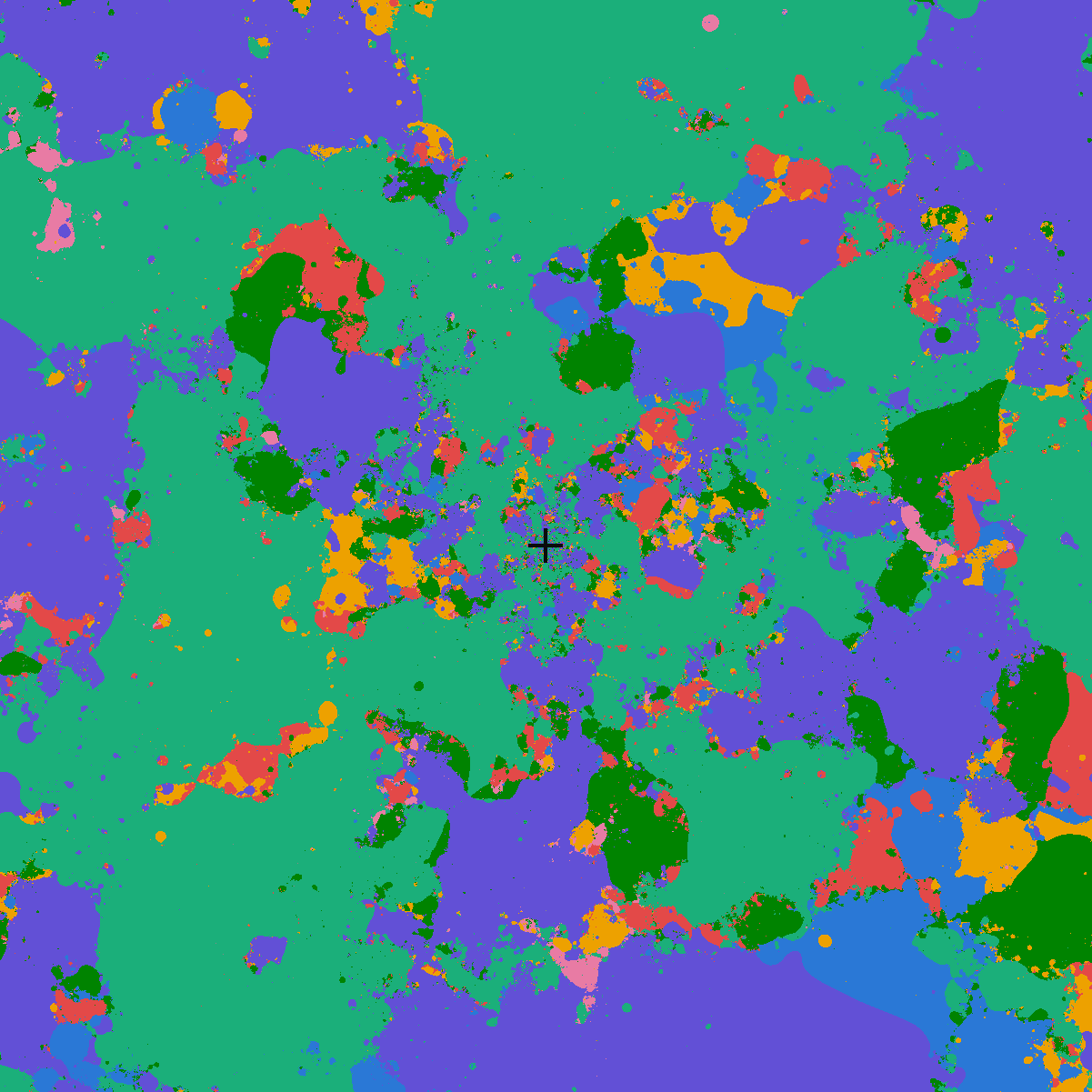}\\[2pt]
\includegraphics[width=0.49\textwidth]{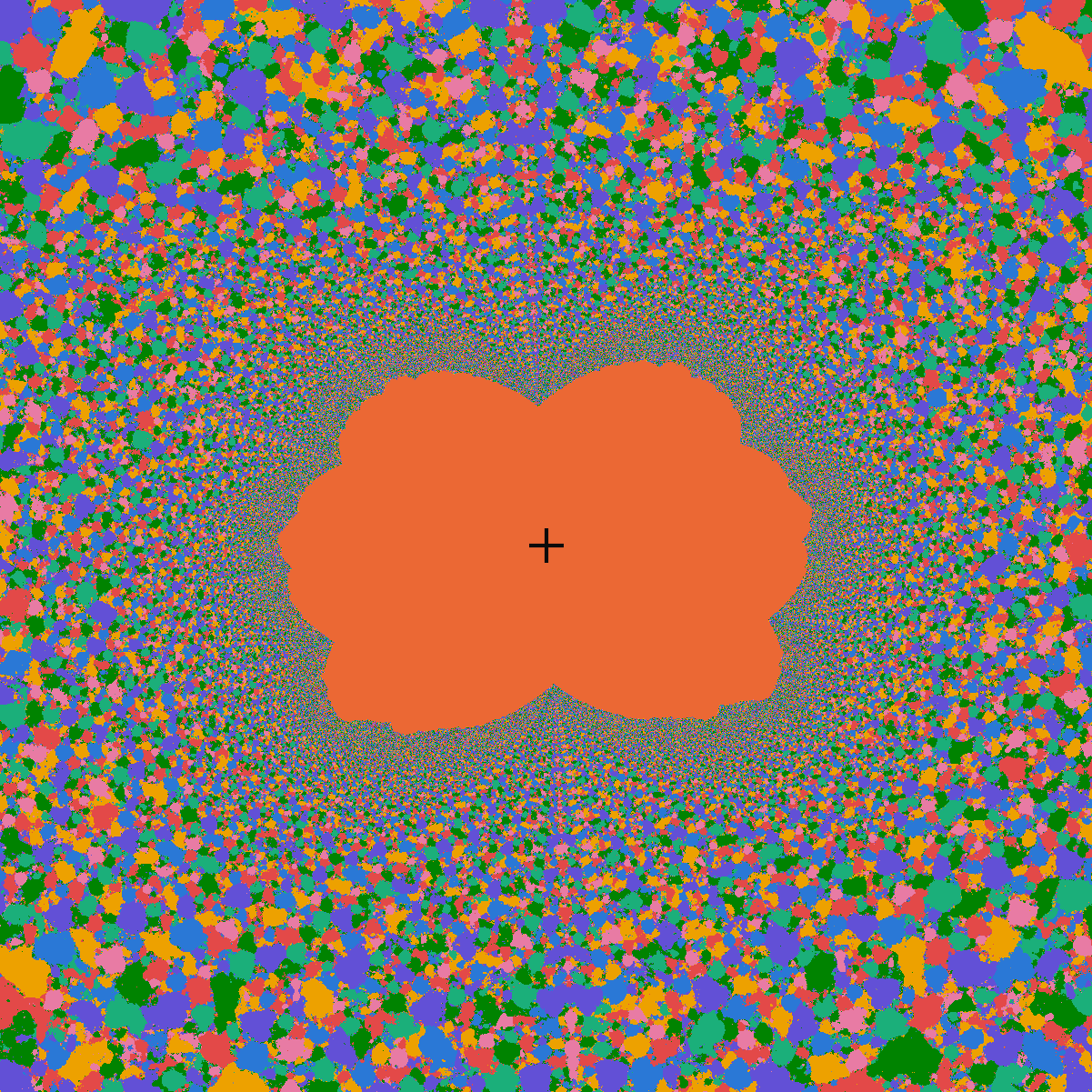}\hspace{2pt}%
\includegraphics[width=0.49\textwidth]{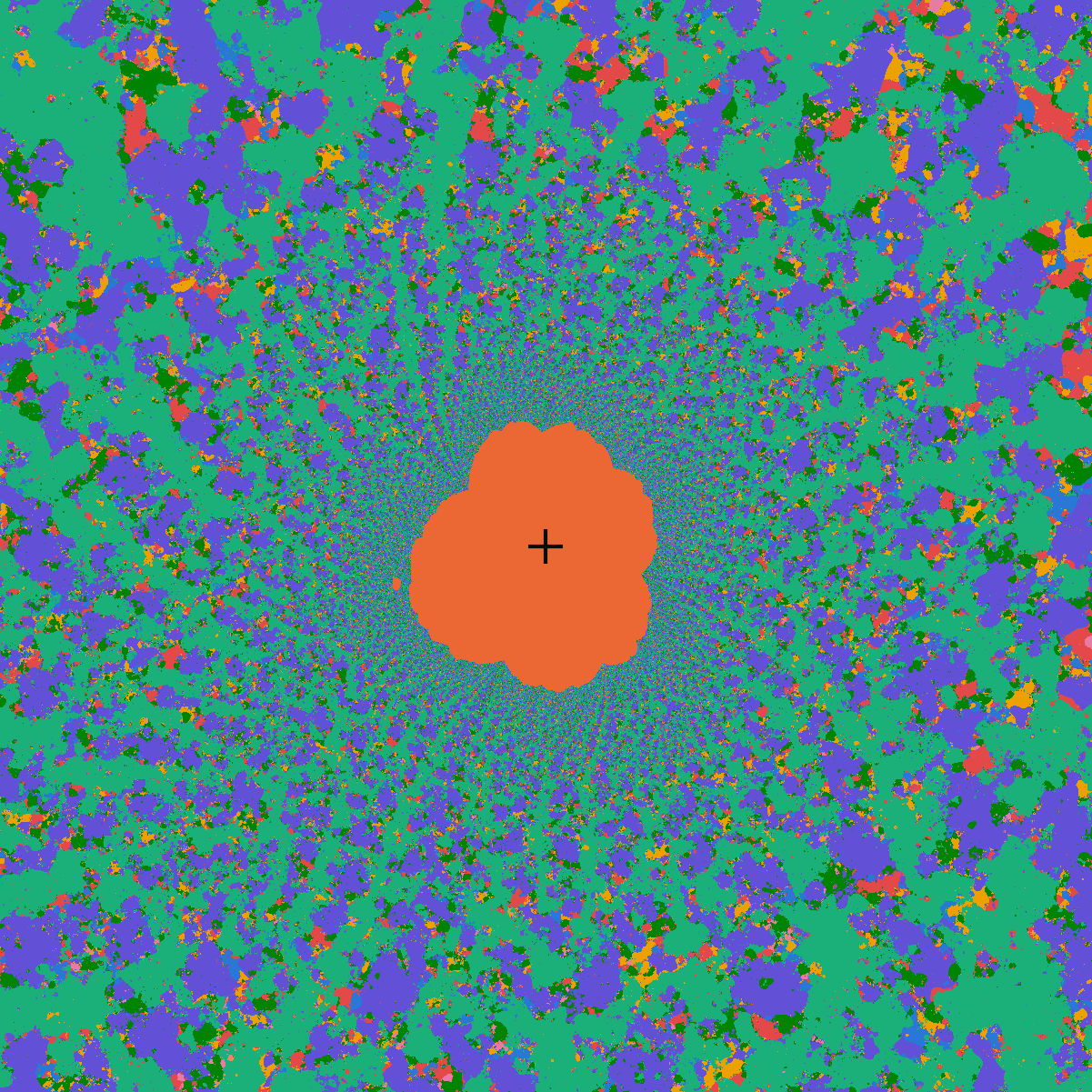}
	\caption{Zooms of the two slices in Figure~\ref{Fig:DynamicSlices1}, at side lengths $2\cdot 10^{-3}$ and $5\cdot 10^{-6}$. Each subsequent zoom is too small to be seen on the previous picture. The basin of the attracting cycle is colored orange in the two pictures at the bottom. The position of the attracting cycle is always indicated by a black cross.}
	\label{Fig:DynamicSlices2}
	\end{figure}

We have tried to measure this quantitatively. For all 14 variables (seven for the positions of the roots, seven for the dynamics), we have plotted a neighborhood in the complex plane to check for which values the attracting cycle persists. Typical plots in the dynamics variables are shown in Figure~\ref{Fig:DynamicSlices1} and Figure~\ref{Fig:DynamicSlices2}, in various magnifications. 

The results of a more quantitative way of measuring are shown in the following table:\\
\medskip
\centerline{
\framebox{
\begin{tabular}{l|rr} 
 & dynamics & roots\\ \hline 
minimal diameter \rule{0pt}{10pt} & $8.93\cdot 10^{-7}$ & $1.855\cdot 10^{-7}$  \\
maximal diameter  & $215.2\cdot 10^{-7}$ & $19.93\cdot 10^{-7}$\\
average diameter  & $30.6\cdot 10^{-7}$  & $5.6\cdot 10^{-7}$ 
\\
relative volume  & $6.3\cdot 10^{-78}$ & $1.02\cdot 10^{-50}$
\end{tabular}
}
}

The left column shows estimates for the diameter of the attracting basin within this coordinate. We have calculated the lengths of the largest horizontal and vertical intervals within the basin (i.e.\ along the real and imaginary directions), and the table presents the maximal and minimal values, as well as the geometric mean and the product of the 14 real lengths. This product is a rough measure of the volume of the attracting basin with respect to a unit cube in $\C^7=\R^{14}$ where the dynamics takes place. The result is extraordinarily tiny: the size of the attracting basin covers only a fraction of about $6.3\cdot 10^{-78}$ of a unit cube in $\C^7$. 

\begin{figure}[htbp]
      \centering
      \includegraphics[width=0.7\linewidth]{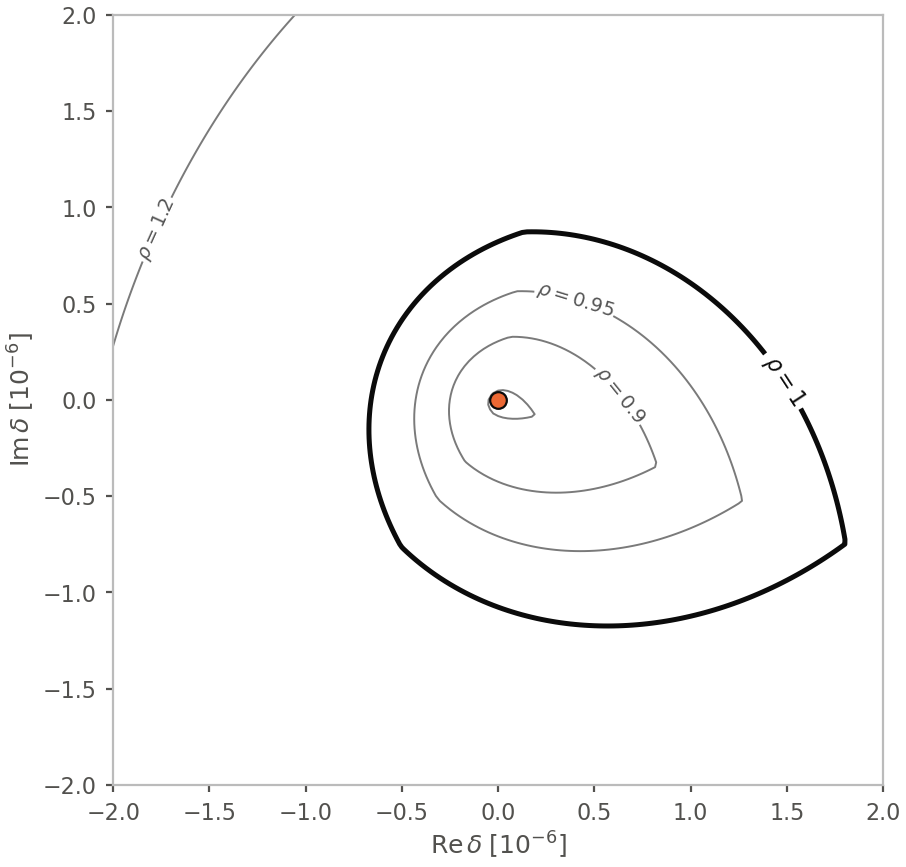}
      \caption{
       A complex one-dimensional slice through the space of polynomials, where the root $r_4$ of $p$ is replaced by $r_4 + \delta$ with the other six roots fixed; shown is the region $|\Re\delta|<2\cdot 10^{-6}$, $|\Im\delta|<2\cdot 10^{-6}$. The picture is centered at our polynomial $p$ with attracting $4$-cycle ($\rho = 0.7353$); the curves show level sets of the largest eigenvalue for the continued cycle. }
      \label{Fig:RootNeighborhood}
  \end{figure}

We have treated the parameter coordinates (i.e.\ the positions of the roots) in a similar way. One difference is that because of the M\"obius equivariance, we only have 4 independent degrees of freedom. Moreover, when changing the positions of the roots, we have followed analytically the periodic orbit that starts out being attracting. A typical picture is shown in Figure~\ref{Fig:RootNeighborhood}. The boundary of the domain where the continued cycle is attracting is indicated (a piecewise analytic curve).

Again, we have measured the diameters of the domain for which an attracting cycle exists, along the $2\times 4=8$ real directions. These are shown in the right column of the table presented earlier. It turns out that the locus in parameter space for which the attracting cycle persists occupies about a fraction of $10^{-50}$ of a unit cube in $\C^4$ that parametrizes the roots.

With respect to the product space of the 7 dynamical and the 4 root coordinates, the locus of attracting cycles occupies thus a fraction of about $6\cdot 10^{-128}$ of a unit cube. In other words, a random point (with respect to Lebesgue measure on $\C^{11}$), when restricted to the ``right'' unit cube, has probability less than $10^{-127}$ to hit a polynomial and a starting vector in the attracting cycle that we found. 

Admittedly the unit cube is not an intrinsically meaningful concept, and all estimates are very rough. In any case, the true search space is certainly much larger than this unit cube, so this explains why an attracting cycle has not been discovered despite the ample successful numerical experiments using the \EA\ method.

%\clearpage

\subsection{Remarks on use of computer clusters and AI}
\label{Sub:RemarksComputers}
The search for attracting cycles was carried out on the Euler Cluster of ETH Z\"urich in eleven computation packages between the 5th and 17th of October 2026, each using $48$ cores for almost $24$ hours, which is about $1.3 \cdot 10^4$ core hours. We thank Kilian Keen for running one of the packages. 
This work was carried out with the help of an AI assistant (Claude, Anthropic), which wrote most of the search and certification code, ran the laptop computations, and took part in analyzing and designing the experiments.

\Newpage

\end{document}